\documentstyle{article}

\newcounter{subequation}
\newenvironment{subequation}%
{\addtocounter{equation}{-1}%
\stepcounter{subequation}%
\begin{equation}}%
{\end{equation}%
}

\newcommand{\beq}{\begin{equation}}
\newcommand{\eeq}{\end{equation}}
\newcommand{\bseq}{\begin{subequation}}
\newcommand{\eseq}{\end{subequation}}
\newcommand{\bea}{\begin{eqnarray}}
\newcommand{\eea}{\end{eqnarray}}
\newcommand{\refeq}[1]{(\ref{#1})}

\newcommand{\noin}{\noindent}

\newcommand{\ol}[1]{\overline{#1}}
\newcommand{\ul}[1]{\underline{#1}}

\newcommand{\pr}{\prime}
\newcommand{\ppr}{{\prime\prime}}

\newcommand{\eps}{\epsilon}

\newcommand{\LL}{\Lambda}

\newcommand{\R}{{\rm I\!R}}
\newcommand{\N}{{\rm I\!N}}
\newcommand{\C}{{\rm I\!\!\!C}}
\newcommand{\Z}{{\rm Z\!\!Z}}
\begin{document}

\title{ CURL-FREE GINZBURG--LANDAU VORTICES}

\author{SAGUN CHANILLO and MICHAEL K.-H. KIESSLING  \smallskip \\
	Department of Mathematics, Rutgers University \smallskip \\
	Piscataway, N.J. 08854\\ \\
	email: chanillo@math.rutgers.edu\\
	email: miki@math.rutgers.edu}

\date{April 23, 1997}

\maketitle

\begin{abstract}
\noindent
For certain nonlinear elliptic PDE problems in two dimensions, 
the classical isoperimetric inequality produces a sharp 
inequality that violates a Pohozaev identity except for 
radial symmetric, decreasing solutions. A generalized
version of this technique is used here to prove 
radial symmetry of curl-free Ginzburg-Landau vortices. 
\end{abstract}

\bigskip
{\bf Keywords:} Ginzburg-Landau equations, radial symmetry,

isoperimetric inequality, Pohozaev identity. 

\bigskip
\bigskip

\centerline{to appear in: NONLINEAR ANALYSIS }

\vfill\eject

\section{Introduction}
Consider a vector field
$u:\R^2\longrightarrow \C$ of class $C^2$ 
that satisfies the system of Ginzburg--Landau equations
\beq\label{GLeq}
-\Delta u(x) = u(1-|u|^2),\ \ \ x\in \R^2
\eeq
together with  the compactness condition
\beq\label{cmpct}
\int_{\bf R^2}(1-|u|^2)^2 dx <\infty \, .
\eeq
Separation of polar coordinates $(r,\theta)$ in the form
\beq\label{svansatz}
u_d(x) = e^{i\, d\theta}f_d(r)\, ,
\eeq
with $d \in \Z$ yields a 
countable set of solutions, with $f_d$ satisfying
\beq\label{fdEQ}
-f_d^\ppr(r) - {1\over r} f_d^\pr + {d^2\over r^2} f_d = 
f_d(1-f_d^2);\ \ r\in (0,\infty)
\eeq
together with the conditions 
\beq\label{fdBC}
f_d(0) = 0 ; \ \ \ \lim_{r\to\infty}f_d(r) = 1\, .
\eeq
Equations \refeq{fdEQ}, \refeq{fdBC} define a unique 
$f_d: \R\mapsto \R^+\cup\{0\}$, 
see \cite{HerveHerve} and also \cite{CohenNeuRosales}, 
\cite{Greenberg}, \cite{Hagan}, \cite{KopellHoward}. 
The symmetry group of the model \refeq{GLeq}, \refeq{cmpct}
acts via rotations and translations on the base space, 
and by complex rotations and conjugations on the target space. 
Thus, every $f_{|d|}(r)$ in fact uniquely defines a whole group 
orbit of solutions. It would be interesting to know
\cite{BethuelBrezisHelein}, 
\cite{BrezisMerleRiviere} whether all solutions of 
\refeq{GLeq}, \refeq{cmpct} are given 
by the group orbits obtained from \refeq{svansatz}.

The problem can be subdivided  according
to the asymptotic winding number, or Brouwer degree, 
 $d = {\rm deg\ } (u,\partial B_R) \in
\Z$ of a solution $u$, which is well
defined for $R$ sufficiently large \cite{BrezisMerleRiviere}.
By a theorem of \cite{BrezisMerleRiviere}, we have
\beq\label{quant}
\int_{\bf R^2}(1-|u|^2)^2 \, dx = 2\pi d^2\, 
\eeq
for any smooth solution of \refeq{GLeq}. It follows from
\refeq{cmpct} and \refeq{quant} that $|d| <\infty$. Furthermore, 
the group action of complex conjugation maps
each solution with winding number $d$ uniquely 
 into one with winding number $-d$. Hence, 
the group orbit of a solution $u$ with winding 
number $d$ contains 
a representative with winding number $|d|$.
Let ${\cal S}$ be the set of group orbits of
solutions of \refeq{GLeq}, \refeq{cmpct}. One
wants to know \cite{BrezisMerleRiviere} 
whether the map ${\cal S}\to \N\cup\{0\}$ is injective.

In \cite{BrezisMerleRiviere} it has been inferred
from Liouville's theorem that the group orbit for
$d=0$ is unique and given by the unimodular constant 
solutions of \refeq{GLeq}, \refeq{cmpct}. 
For $d= 1$, uniqueness of the group orbit
was very recently proved by Mironescu \cite{MironescuCR}.
Moreover, it is known \cite{Shafrir} that all solutions $u$
of \refeq{GLeq}, \refeq{cmpct} 
are asymptotically of the form $u_d$ as given by \refeq{svansatz};
more precisely, by Theorem 1 of Shafrir \cite{Shafrir}, any 
solution $u$ of \refeq{GLeq}, \refeq{cmpct} with degree $d$ satisfies
\beq\label{uasymp}
u(x) - e^{i d\,( \theta - \theta_0)} \to 0 
\eeq
uniformly as $r=|x|\to\infty$, where  $\theta_0$ is 
some constant determined by $u$. Mironescu 
\cite{MironescuCR} announces a theorem 
that degree $d$ solutions that have only one zero belong to
the orbit of \refeq{svansatz}. Whether solutions with several 
zeros exist is still not known.

Our contribution to the characterization of the
solutions of \refeq{GLeq}, \refeq{cmpct} concerns
a  closely related yet different aspect. 
Notice that $u_0$ and $u_{\pm 1}$ are strictly 
radially symmetric in the traditional sense, and
the solutions $u_d$ are radially symmetric  
on the Riemann surface defined by $z \mapsto z^{1/d}$. 
Understood in this way, a ``radially symmetric solution'' 
\cite{JaffeTaubes} is necessarily of the form \refeq{svansatz}. 
The question whether ${\cal S}$ consists of 
the group orbits of \refeq{svansatz} would be 
answered in the affirmative if one can prove 
that each group orbit contains a radially symmetric
solution in the sense of \cite{JaffeTaubes}. 

In \cite{ckGAFA} we proved that certain
conformally invariant systems of nonlinear PDE 
on $\R^2$ can only have radially symmetric solutions. 
Our technique, which is to compare isoperimetric estimates at 
infinity versus Pohozaev's identity, is a generalization
of earlier techniques for a single PDE  in a finite disk 
applied by Bandle and Keady
\cite{BandleBOOK} and by P.L. Lions \cite{Lions};
see also \cite{KesavanPacella} for results in higher dimensions. 
Thus H. Brezis \cite{brezisPC} was led to ask whether 
our technique \cite{ckGAFA} can prove radial symmetry of
the solutions of \refeq{GLeq}, \refeq{cmpct}?

In this paper we introduce an interesting variant of our technique 
\cite{ckGAFA} that does  apply to the class of curl-free
solutions of \refeq{GLeq}, \refeq{cmpct} that have a single zero. 
We prove that any orbit of 
solutions of \refeq{GLeq}, \refeq{cmpct} which contains 
a curl-free vortex with a single zero must be the orbit
of $u_1$, i.e., \refeq{svansatz} with $d = 1$. 
An announcement of this result was made in \cite{ckCR}.
We have not obtained sufficient asymptotic control in $\R^2$ 
of solutions of arbitrary degree $d$ to conclude radial 
symmetry for $|d|>1$ analogs of the curl-free vortices.
This, however, is a technical issue,
not a conceptual one, as will be explained at the end of 
last section. For the related problem with arbitrary finite 
degree $d$ in a disk of radius $R$, 
we prove radial symmetry for those orbits of 
solutions having a single zero at the center of the disk 
which are generated by a curl-free vortex on the Riemann surface
for the map $z\to z^{1/d}$, $z\in B_R\subset \C$. Those 
solutions belong to the analog orbit of \refeq{svansatz} in  $B_R$. 

We emphasize that our technique applies to a wide class of PDE 
problems in two dimensions
and requires only mild control of solutions at infinity. 
In particular, in the problem treated here the a priori divergence 
 at infinity of the solutions is so strong that
we do not see how  to apply the moving plane method \cite{GNN}.

 \vfill\eject

\section{Statement of Results}

We identify the target space $\C$ with $\R^2$. 
For any solution $u$ to \refeq{GLeq}, $SO(2)$ symmetry
implies that ${\cal R}_\beta u$ is also a solution of \refeq{GLeq},
where ${\cal R}_\beta\in SO(2)$ is a rotation by any arbitrary angle 
$\beta \in [0,2\pi)$. By Shafrir's result \refeq{uasymp}, 
if  $u$  has degree $d= 1$ 
we can find a particular $\beta(u)\in [0,2\pi)$ such that
\beq\label{realuasymp}
({\cal R}_{\beta(u)} u)(x)-\frac{x}{|x|} \to 0
\eeq
uniformly in $|x|$. Since 
\beq\label{asymprotfree}
 \nabla\times\left( \frac{x}{|x|} \right) = 0 \, ,
\eeq
we say that our ${\cal R}_{\beta(u)} u$ is asymptotically curl-free. 
While this does not imply that ${\cal R}_{\beta(u)} u$ is globally 
curl-free, it does 
make sense to investigate the subclass of solutions of degree 1 which
are globally curl-free after at most a  rotation. Since this class 
contains the special solution $u_1$ given by \refeq{svansatz},
it is nonempty. The interesting question then is whether 
there exist globally curl-free solutions of degree 1 whose group orbit 
does not contain $f_1$. Our first theorem says that this is not the 
case if the curl-free $u$ has only one zero. Notice below that we 
do not have to impose degree 1.

\medskip
\noin
{\bf Theorem 1. -- }
{\it Let $u(x)$ satisfy \refeq{GLeq}, \refeq{cmpct} and the 
 additional hypotheses that
\beq\label{rotfree}
 \nabla\times({\cal R}_\beta u)=0\ 
\eeq
for some $\beta$, and that
\beq\label{onezero}
 u(x_0)=0 
\eeq
for one and only one $x_0\in \R^2$.
Then 
\beq\label{urepres}
({\cal R}_{\beta}\, u)(x+x_0) =\, \pm {x \over{|x|}} f_1(|x|),
\eeq
with $f_1(|x|)$ being the unique non-negative
solution of \refeq{fdEQ}, \refeq{fdBC}  for $d=1$.}

\medskip
Theorem 1 has been announced in \cite{ckCR}. In section 3 
we give the complete proof. We remark that our assumption 
that $u$ has a single zero is purely technical and needed 
to reduce \refeq{GLeq} globally to a single scalar PDE. 
If one wants to abandon the single zero condition, one also has
to abandon  the idea of reducing \refeq{GLeq} globally
to a scalar PDE. 

For $u(x) = u(z,\ol{z})$ being a solution of \refeq{GLeq}, 
\refeq{cmpct} with $d > 1$, Shafrir's asymptotic 
result \refeq{uasymp} implies that 
$v(x) = u( z^{1/d},\ol{z}^{1/d})$ is well defined for $|z| > R$
and asymptotically curl-free. This observation suggests
that Theorem 1 should have a counterpart for higher degree. 
In section 4, we prove the following related  result in
a disk  $B_R = B_R(0) \subset \R^2$ of radius $R$, 
centered at the origin. 

\medskip
\noin
{\bf Theorem 2. -- }
{\it Let $u(x)$ be a solution of degree $d\in\N$ of
\beq\label{GLeqD}
-\Delta u(x) = u(1-|u|^2),\ \ \ x\in B_R(0),
\eeq
with
\beq\label{Dbc}
u(x) = e^{id\,\theta}, \ \ \ x \in \partial B_R(0)\, ,
\eeq
satisfying the additional hypotheses that
\beq\label{onenull}
 u(x_0)=0 
\eeq
if and only if $x_0=0$, and that
\beq\label{curlfree}
 \nabla\times v(x)=0\, , 
\eeq
where $v(x)=v(z,\ol{z})$ is given by 
\beq
v(z,\ol{z}) = u \bigl( z^{1/d},\, \ol{z}^{1/d}\bigr) \, .
\eeq
Then 
\beq\label{uSOL}
u(x) =\, e^{id\,\theta} F_d(|x|)\, ,
\eeq
where $F_d(|x|)$ is the unique non-negative solution of 
\beq\label{FdEQ}
-F_d^\ppr(r) - {1\over r} F_d^\pr + {d^2\over r^2} F_d = 
F_d(1-F_d^2);\ \ r\in (0,1)
\eeq
satisfying
\beq\label{FdBC}
F_d(0) = 0 ; \ \ \ F_d(R) = 1\, .
\eeq
}
\medskip

\noin
{\bf Remark:} When $d= 1$  the condition $x_0 =0$ can be dropped.

We conclude with comments on the analogous problem in $\R^2$.

\section{Proof of Theorem 1} 

To begin the proof of Theorem 1, we state two simple lemmata.

\medskip
\noin
{\bf Lemma 1. -- } 
{\it Let $u$ satisfy \refeq{GLeq}, \refeq{cmpct}, and 
the hypotheses \refeq{rotfree} of Theorem 1. 
Then $u$ is given by 
\beq\label{gradrep}
{\cal R}_{\beta} u = \mp \nabla\phi,
\eeq
for some scalar-valued function $\phi(x)$ which satisfies
\beq\label{phiASYMP}
\lim_{|x| \to \infty} |x|^{-1} \phi(x) = -1
\eeq
uniformly as $|x| \to \infty$. }
\medskip

{\it Proof of Lemma 1: } By hypothesis, \refeq{rotfree} holds for some
$\beta$. Then there exists a scalar function $\phi$ such that
${\cal R}_{\beta} u = \mp \nabla\phi$. 
It further follows from  the asymptotics \refeq{uasymp} 
of Shafrir \cite{Shafrir} that
$|\phi(x)|=\,|x|+o(|x|)$, uniformly as $|x|\to\infty$. 
By $SO(2)$ invariance of \refeq{GLeq}, we are free to choose
\beq\label{phiasymp}
\phi(x) =\, -|x|+o(|x|)
\eeq
uniformly as $|x|\to\infty$.  This proves \refeq{phiASYMP}. QED
\bigskip

\medskip
\noin
{\bf Lemma 2. -- }{\it Let $u$ satisfy the hypotheses  of Lemma 1. 
such that $u$ is represented by \refeq{gradrep}. Then $|\nabla\phi|$ is
constant on the level curves of $\phi$. }
\medskip

{\it Proof of Lemma 2: } Since ${\cal R}_\beta u$ is also a solution
to \refeq{GLeq}, we may substitute the right side of 
\refeq{gradrep} into \refeq{GLeq} to get,
\beq\label{GLphi}
  \nabla \big(-\Delta \phi (x)\big) = \, \lambda(x)\nabla\phi(x),
\eeq
where, 
\beq\label{lambdaofphi}
\lambda(x)=\, 1-|\nabla\phi(x)|^2.
\eeq
We now take the curl of \refeq{GLphi} and, since $\nabla\times\nabla f =0$ 
for  a scalar function $f$, we get
\beq\label{lambdatimesphi}
\nabla \lambda(x) \times \nabla\phi(x) = 0.
\eeq
By \refeq{lambdatimesphi}
and \refeq{onezero} it follows that $\lambda(x)$ is constant on 
the level curves of $\phi$. QED
\bigskip

The single zero hypothesis enters next  and 
leads to the conclusion that $\lambda(x)$ 
is in fact a function of $\phi(x)$. As a consequence,
for solutions of \refeq{GLeq} of the type \refeq{rotfree},
we can reduce \refeq{GLeq} to a single scalar equation 
valid globally in $\R^2$.

\medskip
\noin
{\bf Lemma 3. -- }{\it Let $u$ satisfy the hypotheses of Theorem 1,
such that $u$ is represented by \refeq{gradrep}. Then,
globally in $\R^2$, the function $\phi$ solves the scalar equation
\beq\label{phiPDE}
-\Delta \phi = g(\phi)\ 
\eeq
for some increasing differentiable function
 $g: \R\to\R^+$ that satisfies
\beq\label{gasymp}
\lim_{t\to -\infty} t g(t) = -1
\eeq
}
\medskip

\noindent
{\bf Remark:} Lemma 3 does not imply that there is a universal
function $g$ for all $u$ that satisfy the hypotheses of Theorem 1.
Each such $u$ may lead to its own $g$. 
\medskip

{\it Proof of Lemma 3: }  We recall that $u$ is a smooth vector field
with \cite{BrezisMerleRiviere} 
\beq\label{absubound}
|u| \leq 1 \ 
\eeq
everywhere in $\R^2$, with
\beq\label{absuto1}
|u| \to 1 
\eeq
uniformly as $|x|\to\infty$, and with
\beq\label{absnabubound}
|\nabla u| \in L^\infty
\eeq
at all points of $\R^2$. Thus, $|\nabla\phi|\to 1$ 
uniformly as $|x|\to\infty$. Whence, $\phi$ has 
critical points only for $|x|<\, R$. 
By hypothesis \refeq{onezero}, $\phi$
has only one critical point in all of $\R^2$. 
By \refeq{phiASYMP}, each level curve $\{x:\phi(x)=\, t\}$ is compact. 
The facts that $\phi$ has only one critical point, which is not
at infinity, and that level curves are compact 
ensures that for a regular value $t$ for $\phi$ 
each level curve $\{x:\phi \, =t\}$ has at most one connected 
component and no self-intersections. In addition, by 
\refeq{onezero} the set of critical values for $\phi$ 
consists of a single point. By \refeq{onezero} again, this 
critical value is attained at a single point of $\R^2$.
Now denote by ${\cal S}$ the range of $\phi$. ${\cal S}$ is
clearly a sub-interval of $\R$. The discussion above shows
that we have a well-defined function $\Theta:{\cal S}\to\R$, 
such that $\lambda(x)= (\Theta\circ\phi)(x)$, i.e.,  with \refeq{lambdaofphi},
\beq\label{lambdaTheta}
1- |\nabla \phi(x)|^2 = (\Theta\circ\phi)(x) \, .
\eeq
Moreover, \refeq{absubound} implies $|\nabla \phi| \leq 1$, which 
clearly implies $\Theta \geq 0$. 

The discussion above renders \refeq{GLphi} in the form
\beq\label{nablaphieq}
\nabla \big(-\Delta \phi (x)\big) =\, (\Theta\circ\phi)(x)\nabla\phi(x)\, .
\eeq
We now see right away 
from \refeq{nablaphieq} that we can find a differentiable 
function $g:\R\to\R$, with $g^\pr(t)=\, \Theta(t)$, such
that $\phi$ satisfies the scalar equation \refeq{phiPDE} in $\R^2$. 
In fact, we can choose
\beq\label{gdefine}
g(t) = \int\limits_{-\infty}^t \Theta(s) ds\, .
\eeq
The  result of  Shafrir  \cite{Shafrir} that
\beq\label{moduasymp}
1- |u(x)|^2 = {d^2\over r^2} + o\left({1\over r^2}\right)\, ,
\eeq
uniformly as $|x|\to\infty$, 
shows that the integral  \refeq{gdefine} exists. 
Moreover, since $\Theta \geq 0$, we conclude that $g$ is
nonnegative and increasing. Finally, \refeq{moduasymp}
implies \refeq{gasymp}. This completes the proof of Lemma 3. QED
\bigskip

We will now prove the radial symmetry of solutions of \refeq{phiPDE}, 
\refeq{phiASYMP}, \refeq{lambdaTheta} by inferring from 
the isoperimetric inequality that any hypothetical
nonradial solution would violate an a-priori identity that is 
satisfied by all solutions. In \cite{ckGAFA} we could use
Pohozaev's identity. Here we have to take a more complicated 
identity based  on Pohozaev's. This 
interesting variant of the strategy of \cite{ckGAFA}
is required by the divergence \refeq{phiasymp} of 
$\phi$ at infinity, which is too strong for the 
arguments in \cite{ckGAFA} to apply directly.

For each solution $\phi(x)$ of \refeq{phiPDE}, \refeq{phiASYMP}, 
we  define its non-increasing radial re-arrangement,
centered at the origin, denoted by $\phi^*(|x|)$.
Let $\Lambda_t = \{x: \phi (x) > t \}$, 
so that $\partial\Lambda_t = \{x:\phi(x)=\, t\}$. 
Since, by \refeq{phiASYMP} and \refeq{onezero}, each level curve 
$\partial\Lambda_t$ is compact and smooth, 
the map $t\mapsto |\Lambda_t|$ is well defined and decreasing.
It diverges to $\infty$ for $t\to -\infty$. Let 
$B_r = \{x: |x|<r \}$, and define $\rho(t)$ by 
\beq\label{raddistrfunc}
|B_{\rho(t)}| = |\Lambda_t| \, .
\eeq  
Then $\phi^*(r)$ is defined as the  radial symmetric non-increasing 
function that satisfies $B_{\rho(t)} = \{x: \phi^* (\rho(t)) > t \}$.

With the help of $\phi^*(|x|)$, 
we uniquely assign to each $\phi(x)$ the functions 
\beq\label{Mdefin}
M(r)=\, 2\pi\int_0^r (g\circ\phi^*)(s)s\, ds, 
\eeq
\beq\label{Adefin}
A(r)=\, 2\pi\int_0^r (G\circ \phi^*)(s)s\, ds,
\eeq
where  $G$ is a primitive of $g$, i.e., $G^\pr(t)=g(t)$, and finally
\beq\label{Hdefin}
H(r)= \, {1\over 2}\big(M(r)\big)^2 +
 2\pi r^3 {d \over d r}\biggl({{A(r)}\over{r^2}}\biggr)\, .
\eeq
Note that $H(r)$ is independent of the choice of the primitive $G$.

We are now in a position to state 
two key lemmata about solutions of \refeq{phiPDE}.
The first one states an inequality which is a consequence
of  the classical iso-perimetric inequality
and the Cauchy-Schwartz inequality. 

\medskip
\noin
{\bf Lemma 4. -- } 
{\it Suppose $\phi$ satisfies \refeq{phiPDE} and \refeq{phiASYMP}. Then 
\beq
\lim_{r\to\infty}H(r)\geq 0\, ,
\eeq
with equality holding if and only if $\phi$ is radial symmetric and 
decreasing about some point.}
\medskip

On the other hand, with the help of Pohozaev's identity we find 
the following identity valid for {\bf all} solutions $\phi$ of 
\refeq{phiPDE}, \refeq{phiASYMP}, \refeq{lambdaTheta}.

\medskip
\noin
{\bf Lemma 5. -- } 
{\it Suppose $\phi$ satisfies \refeq{phiPDE} and \refeq{phiASYMP},
with $|\nabla \phi|$ satisfying \refeq{lambdaTheta}.  Then we have
\beq
\lim_{r\to\infty}H(r) = 0.
\eeq
}
\medskip

 Lemmata 4 and 5 immediately imply the following. 

\medskip
\noin
{\bf Proposition 1. -- } 
{\it Suppose $\phi$ satisfies equation \refeq{phiPDE} in $\R^2$,
with asymptotics given in \refeq{phiASYMP}, and
with $|\nabla \phi|$ satisfying \refeq{lambdaTheta}. 
Then there exists some center of symmetry $x_0$ such that
\beq
\phi(x)= \phi^*(|x-x_0|)\, .
\eeq
}
\medskip

Now notice that if ${\cal R}_\beta u=\mp \nabla\phi$ 
satisfies \refeq{GLeq} and \refeq{cmpct}, with 
$\phi$ satisfying \refeq{phiPDE} and \refeq{phiASYMP},
then Proposition 1 applies to $\phi$ as 
a result of Lemmata 1 - 5. By translation invariance,
the center of symmetry $x_0$ can be identified with the
origin. Then $\phi(x) = \phi^*(r)$, and ${\phi^*}^\pr(r)$ satisfies
\refeq{fdEQ}, \refeq{fdBC}, which in turn has a unique
solution, by \cite{HerveHerve}. Our Theorem 1 follows. 
\medskip

It remains to present the proofs of Lemma 4 and Lemma 5.

\medskip
{\it Proof of Lemma 4. } We begin with the known result,
see \cite{BandleBOOK}, that
\beq\label{isopestim}
- 2\pi r {\phi^*}^\prime (r)\Bigl\vert_{r = \rho(t)} \Bigr.
\leq - \int_{\partial\Lambda_t}
\langle\nu,\nabla \phi \rangle d\sigma
\eeq
with equality holding if and only if 
$\phi^*(|x-x_0|) = \phi(x)$ for some $x_0$.
Applying Green's theorem and then equation \refeq{phiPDE}
to the r.h.s. of \refeq{isopestim} we get
\beq\label{greenPDE}
- \int_{\partial\Lambda_t} \langle\nu,\nabla \phi \rangle d\sigma =
\int_{\Lambda_t}g( \phi) dx \, .
\eeq
By equi-measurability of $\phi$ and $\phi^*$, 
and noting \refeq{Mdefin}, we have
\beq\label{equimeasg}
\int_{\Lambda_t}g( \phi) dx = \int_{B_{\rho(t)}} g(\phi^*) dx =
M(\rho(t)) \, .
\eeq
Combining \refeq{isopestim}, \refeq{greenPDE} and \refeq{equimeasg} gives
\beq\label{isor}
- 2\pi r {\phi^*}^\prime (r) \leq M(r) \, .
\eeq
Now
\beq\label{Aprpr}
A^{\prime\prime}(r) = 
2\pi [r g(\phi^*) {\phi^*}^\prime(r) + G(\phi^*)]\, .
\eeq
Inserting \refeq{isor} into \refeq{Aprpr}, 
and recalling $g>0$, we  get
\beq\label{Aprprestim}
2\pi r A^{\prime\prime} (r) - 2\pi A^\prime (r)
+ M(r) M^\prime(r)  \geq 0
\eeq
for all $r$. By using 
\beq\label{diffidentity}
2\pi r A^{\prime\prime} (r) - 2\pi A^\prime(r)  = 
{d\over d r}\left(2\pi r^3 {d\over d r}\left(\frac{A(r)}{r^2}\right) \right)
\eeq 
and comparing with the derivative of \refeq{Hdefin}, we find
\beq\label{hprimpos}
H^\prime(r)   \geq 0\, .
\eeq 
From the definition \refeq{Hdefin} of $H$ we also 
easily find $H(0) = 0$. Hence, upon integration of \refeq{hprimpos},
\beq
\liminf_{r\to\infty} H(r) \geq 0 \, ,
\eeq
with equality holding if and only if 
$\phi^*(|x-x_0|) = \phi(x)$ for some $x_0$. QED
\bigskip

{\it Proof of Lemma 5. } In \refeq{Hdefin} we assigned 
a unique function $H(r)$  to each solution $\phi(x)$
by integrating over its radial rearrangement $\phi^*(|x|)$.
We can also think of $H(r)$ as a function of the level value 
$t$ of $\phi(x)$ on $\partial\Lambda_t$. Formally, this is
achieved by setting $r = \rho(t)$ in $H(r)$, where $\rho(t)$ 
is the unique radius assigned by \refeq{raddistrfunc}
to each level curve $\partial\Lambda_t =\{x:\phi(x) = t\}$. 
To obtain a more transparent expression for $H(\rho(t))$, 
we rewrite the two terms in the r.h.s. of the definition 
\refeq{Hdefin} as integrals over $\partial \Lambda_t$.

Recall the standard form of
the Rellich--Pohozaev identity for \refeq{phiPDE} on $\Lambda_t$, 
\beq\label{poho}
{1\over 2}
\int_{\partial \Lambda_t} \langle x,\nu\rangle |\nabla \phi|^2 d\sigma 
+ \int_{\partial \Lambda_t} \langle x,\nu\rangle G(\phi) d\sigma -
2\int_{\Lambda_t} G( \phi)  dx  = 0\, ,
\eeq
which is obtained as usual by integrating 
$- \langle x, \nabla\phi\rangle \Delta\phi$ over $\Lambda_t$. 
We rewrite the second integral in \refeq{poho}, using
$\phi |_{\partial \Lambda_t} = t$ to  pull out $G(t)$ from the integral,
then using Green's theorem and $\nabla\cdot x\,=2$ to obtain
\beq\label{secinta}
\int_{\partial \Lambda_t} \langle x,\nu\rangle G(\phi) d\sigma =
2 G(t) |\Lambda_t| \, .
\eeq
But 
\beq\label{secintb}
2 G(t) |\Lambda_t|  = 
2 \bigl(G(\phi^*(r)) |B_r|\bigr)\Bigl\vert_{r=\rho(t)}\Bigr.
 = \rho(t) A^\prime (\rho(t)) \, ,
\eeq
which obtains by differentiating \refeq{Adefin} w.r.t. $r$, then using 
$|\Lambda_t| = |B_{\rho(t)}|$. We also rewrite 
the third integral in \refeq{poho}, using equi-measurability of 
$\phi$ and $\phi^*$, and then \refeq{Adefin}, to get
\beq\label{thiint}
2 \int_{\Lambda_t} G( \phi)  dx =
 2\int_{B_{\rho(t)}} G( \phi^*)  dx = 2 A({\rho(t)}) \, .
\eeq
Using now
\beq
r A^\prime(r) - 2 A(r) = r^3 {d\over d r}\left( {A(r)\over r^2}\right) \, ,
\eeq
we see that we can rewrite \refeq{poho} as
\beq\label{pohorew}
\left(r^{3} {d\over d r}\left( {A(r)\over r^2} 
\right)\right)\Bigl\vert_{r=\rho(t)}\Bigr. = 
- {1\over 2}
\int_{\partial \Lambda_t} \langle x,\nu\rangle |\nabla \phi|^2 d\sigma
\, .
\eeq

On the other hand, reading \refeq{equimeasg} and \refeq{greenPDE} backwards, 
and recalling $g >0$, we obtain 
\beq\label{Mofrhot}
M(\rho(t)) 
=   \int_{\partial \Lambda_t} |\nabla \phi| d\sigma \, . 
\eeq
Squaring \refeq{Mofrhot}, then multiplying by 1/2 and
adding $2\pi$ times \refeq{pohorew} gives 
\beq\label{Hsurfint}
 H(\rho(t)) = 
{1 \over 2} \left(\int_{\partial \Lambda_t} |\nabla \phi| d\sigma\right)^2 
-\pi \int_{\partial \Lambda_t} \langle x,\nu\rangle |\nabla \phi|^2 d\sigma
\, .
\eeq

We now rewrite \refeq{Hsurfint} further by making convenient use
of \refeq{lambdaTheta} and \refeq{gdefine}, which state that 
$|\nabla \phi|^2 = 1 - \Theta(t)$ on $\partial\Lambda_t$.
Pulling out the  constant $|\nabla\phi|$ 
terms from both integrals in \refeq{Hsurfint},
then using Green's theorem and $\nabla\cdot x\,=2$ to get
\beq\label{twoarea}
\int_{\partial \Lambda_t} \langle x,\nu\rangle d\sigma = 
\int_{\Lambda_t} \nabla\cdot x\, dx =  2 |\Lambda_t| \, ,
\eeq
we obtain
\beq\label{Hrew}
H(\rho(t)) =  {1 \over 2}  (1-\Theta(t)) 
\left(|\partial\Lambda_t|^2  - 4\pi | \Lambda_t | \right)\, .
\eeq
By \refeq{raddistrfunc} we have $|\Lambda_t| =  \pi \rho^2(t)$.
This implies 
\beq\label{areadiff}
\frac{d}{dt}|\Lambda_t| =  2\pi \rho(t)\rho^\pr(t)\, .
\eeq
On the other hand, the co-area formula and the constancy of  $|\nabla\phi|$ on
$\partial\LL_t$ give us
\beq\label{coareadiff}
\frac{d}{dt}|\Lambda_t| = 
- \int_{\partial\LL_t} {d\sigma_t \over |\nabla\phi|}
= - {|\partial\LL_t| \over \sqrt{1-\Theta(t)}}\, .
\eeq
Hence, \refeq{Hrew} becomes
\beq\label{Hrewr}
H(\rho(t)) =  2\pi^2 \bigl(1-\Theta(t)\bigr) \rho(t)^2
\Bigl((1-\Theta(t))\rho^\pr(t)^2 - 1\Bigr)\, .
\eeq

Now notice that as $t\to-\infty$,  
by \refeq{gasymp} and \refeq{gdefine} we have
$1-\Theta(t)\sim 1 - t^{-2} + o(t^{-2})$, 
and by \refeq{phiASYMP} we have $\rho(t) \sim - t + o(t)$. 
Therefore, to prove Lemma 5, we need to show that
the term exhibited in large parentheses in \refeq{Hrewr} 
is $o(t^{-2})$. 

Now suppose not, that the last term in \refeq{Hrewr} is
$o(t^{-2})$. Then there exists an $a\neq 0$ such that
\beq\label{ODinEsq}
(1-\Theta(t))\rho^\pr(t)^2 \geq 1 + a^2 t^{-2}\, .
\eeq
Using $1-\Theta(t)\sim 1 - t^{-2} + o(t^{-2})$ 
and $\rho^\pr(t)= -1/|\nabla\phi^*(\rho(t))| <0$, 
we infer from \refeq{ODinEsq} that
\beq\label{phiSTARinE} 
-\rho^\pr(t) = 1/|\nabla\phi^*(\rho(t))|
\geq 1 + \frac{1+a^2}{2t^2} + o(t^{-2}) \, .
\eeq

Now let $\ol{r}(t)$ be the radius of the smallest disk
centered at $x_0$ containing $\LL_t$, $\ul{r}(t)$ the radius of the 
largest disk that is centered at $x_0$ and contained in $\partial\LL_t$. 
The locations of the points at which
the outer and inner disks touch $\partial\LL_t$, 
are denoted by $\ol{x}(t)$ and $\ul{x}(t)$, respectively. 
Equi-measurability of $\LL_t$ 
and $B_{\rho(t)}$ implies the ordering
$\ol{r}(t) \geq \rho(t) \geq \ul{r}(t)$.
Therefore, for all $t$, 
\beq\label{distances}
\ol{r}(t) - \ul{r}(t) \geq \rho(t) - \ul{r}(t)\, .
\eeq
We now look at the derivatives of the left and right sides of
\refeq{distances} for large negative $t$, using the fact
that locally $\ol{x}(t)$ and $\ul{x}(t)$ 
are transported along an integral curve of $\nabla\phi$. 
Therefore, and by \refeq{moduasymp}, we have
\beq\label{deriulr}
\frac{d}{dt}\ul{r}(t) =  - \frac{1}{|\nabla\phi|(\ul{x}(t))} 
= - 1 -  \frac{1}{2\ul{r}(t)^2} + o\left(\frac{1}{\ul{r}^2}\right)\, . 
\eeq
An analogous formula holds for the derivative of $\ol{r}(t)$. 
By \refeq{phiasymp}, $\phi(x)= - |x| + o(|x|)$,
so that in leading order we can replace $\ol{r}(t)^2$ 
and $\ul{r}(t)^2$ by $t^2$. Hence, 
\beq\label{deriA}
\frac{d}{dt}\left(\ol{r} - \ul{r}\right)(t) = 
o\left(\frac{1}{{t}^2}\right) \, . 
\eeq
This means, for any small $\epsilon$, 
we can find a large negative $\tau$ such that,
for $t<\tau$ we have
\beq\label{intderiA}
\left(\ol{r} - \ul{r}\right)(t)= 
\int_{-\infty}^t  o\left(\frac{1}{{\xi}^2}\right)d\xi 
\leq \frac{\epsilon}{|t|}
\eeq
In particular, we can choose $\epsilon = 10^{-10}a^2$. 
On the other hand, by \refeq{deriulr} and
\refeq{phiSTARinE} we have
\beq\label{deriB}
\frac{d}{dt}\left({\rho} - \ul{r}\right)(t) \geq 
\frac{a^2}{2t^2} + o(t^{-2})\, ,
\eeq
which upon integration from $-\infty$ to $t$
gives, for all $t< \tau$, 
\beq\label{intderiB}
\left({\rho} - \ul{r}\right)(t) \geq
\int_{-\infty}^t \frac{a^2}{{\xi}^2} + o\left(\frac{1}{{\xi}^2}\right)d\xi 
= \frac{a^2}{|t|} + o(t^{-2})\, .
\eeq
Clearly, \refeq{intderiB} and \refeq{intderiA} lead to the
inequality
\beq\label{reversedistances}
\lim_{t\to -\infty} \frac{\ol{r}(t) - \ul{r}(t)}{\rho(t) - \ul{r}(t)} 
\leq \frac{\epsilon}{a^2} = 10^{-10} \, .
\eeq
But \refeq{reversedistances} contradicts \refeq{distances}.
Hence,  $a = 0$, and it follows that
the term exhibited in large parentheses in \refeq{Hrewr} 
is $o(t^{-2})$. Thus, 
\beq\label{Hatinfty}
\lim_{t\to - \infty}  H(\rho(t)) = 0\, ,
\eeq
and Lemma 5 is proved. QED
\bigskip


\subsection{Comments on the Proof}

In the proof of Lemma 5, we made convenient use of 
\refeq{lambdaTheta} to rewrite \refeq{Hsurfint} into \refeq{Hrewr}.
However, we wish to emphasize that our general strategy of proving 
radial symmetry is not based on \refeq{lambdaTheta} in any 
essential manner. It operates with \refeq{Hsurfint} and is
designed for the general case of PDE of the type \refeq{phiPDE}
with asymptotically radial data for $\phi$, which may be diverging
to $-\infty$. Global radial symmetry, and decrease, of $\phi$ follows
whenever one can show that \refeq{Hsurfint} tends $\to 0$ 
as $\rho(t)\to \infty$. In particular,
at the end of the next section we
raise an open question on the asymptotic control of
the degree $d$ generalization of curl-free vortices. 
An affirmative answer would entail 
that $\liminf_{t\to\infty} H(\rho(t)) \leq 0$, and 
radial symmetry would follow once again by comparing with
Lemma 4. 

On the other hand, the constancy of $|\nabla\phi|$ 
on $\partial\Lambda_t$, here a byproduct of the 
Ginzburg--Landau equations, allows us to 
give a different, more direct proof of Theorem 1.
For the sake of completeness, we include it here.
This proof can be generalized to higher dimensional 
problems of the same type. 

\medskip
\noin
{\it Alternate Proof of Theorem 1:}
We study the integral curves of the vector field $\nabla\phi$ that leave
$x_0$, the critical point of $\phi$, and hit any point $y$ 
on $\partial\LL_t$. By the constancy of $|\nabla\phi|$ on 
the level curves, the length of such an integral curve is 
a function only of $t$, independent of $y$, given by 
\beq\label{curvelength}
L_{x_0,\partial\LL_t}(\nabla\phi) = 
\int_t^{t_0} \frac{d \xi}{\sqrt{1-\Theta(\xi)}}
\eeq
where $t_0 = \phi(x_0)$. Here we used that
$|\nabla\phi|^2_{\partial\LL_t} =  1-\Theta(t)$. 
As before, let $\ol{x}(t)$ be the point on $\partial\LL_t$ 
that has maximum distance from $x_0$, and denote
that distance by $\ol{r}(t)$. 
Since $L_{x_0,\partial\LL_t}(\nabla\phi)$ is independent of
the point $y\in \partial\LL_t$ through which the integral curve
of $\nabla\phi$ runs, we can choose $y=\ol{x}(t)$.
Obviously we have
\beq\label{distanceestim}
 L_{x_0,\partial\LL_t}(\nabla\phi) \geq \ol{r}(t)\, .
\eeq

We now estimate $1-\Theta(t)$ in terms of $|{\phi^*}^\pr|$. 
By \refeq{raddistrfunc}, the right sides of 
 \refeq{coareadiff} and of
\beq\label{coareadiffstar}
\frac{d}{dt} |B_{\rho(t)}| = 
- \int_{\partial B_{\rho(t)}} {d\sigma_t \over |{\phi^*}^\pr|}
\eeq
are equal. Since $\nabla\phi$ has only a 
single zero, by \refeq{onezero}, and this zero is at $x_0$, we 
have that $|\nabla\phi|_{\partial\LL_t} \neq 0$ away from $x_0$,
and also $|\nabla\phi^*|_{\partial B_{\rho(t)}} \neq 0$ away from
$x_0$. Using this,
and again the facts that $|\nabla\phi|$ is constant on $\partial\Lambda_t$ 
and $|{\phi^*}^\pr(\rho(t))|$ on $\partial B_{\rho(t)}$, away from $x_0$
we find from
\refeq{coareadiff}, \refeq{coareadiffstar} and \refeq{raddistrfunc}
that
\beq\label{coareadiffequiv}
{ |\nabla{\phi^*}|_{\partial B_{\rho(t)}}}{|\partial\Lambda_t|} =
{|\nabla\phi|_{\partial \Lambda_t}{|\partial B_{\rho(t)}| }}\, .
\eeq
The classical isoperimetric inequality in its weak form
simply states that $|{\partial\LL_t}| \geq |{\partial B_{\rho(t)}}|$.
Applying this to \refeq{coareadiffequiv} gives
\beq\label{gradcompare}
|\nabla\phi|_{\partial\LL_t} \geq |\nabla\phi^*|_{\partial B_{\rho(t)}} \, .
\eeq
Using \refeq{gradcompare} in \refeq{curvelength}, we find 
the upper bound on $L_{x_0,\partial\LL_t}(\nabla\phi)$ given by
\beq\label{lengthestim}
L_{x_0,\partial\LL_t}(\nabla\phi) 
\leq\int_t^{t_0} \frac{d
\xi}{|{\phi^*}^\prime(\rho(\xi))|} = \rho(t)\, .
\eeq

By \refeq{distanceestim} and \refeq{lengthestim}, we infer that
$\LL_t \subset B_{\rho(t)}(x_0)$ for all $t$. By equimeasurability, 
this now means $\LL_t = B_{\rho(t)}(x_0)$ for all $t$. Therefore,
the level curves of $\phi$ are concentric circles. From here
the rest of Theorem 1 follows immediately. QED. 
\medskip

The above proof, which is based on the constancy of 
$|\nabla\phi|$  on the level curves in an essential way, 
still uses the single zero hypothesis. 
J.J. Aly has suggested to us that by 
methods due to Pucci and Serrin it may be
possible to drop the single zero hypothesis \refeq{onezero}. 
We mention here that dropping the single zero hypothesis 
implies, however, abandoning the single scalar PDE \refeq{phiPDE}, too. 
This raises interesting new problems. 
We hope to address them at a later time. 

\section{Proof of Theorem 2}

To map $u$ to the Riemann surface indicated in Theorem 2, we need the
following lemma. 

\medskip
\noin
{\bf Lemma 6. -- } 
{\it Let $\Omega\subset \C$ be a domain containing the origin.
Assume $u: \Omega \to \C$ is a solution of degree
$d\in\N\cup\{0\}$ of \refeq{GLeqD}, \refeq{onenull}. Then
\beq
u (z,\ol{z}) = a_1 z^d + a_2 \ol{z}^d + O\left(|z|^{d+1}\right)
\eeq
where either $a_1\neq 0$ or $a_2\neq 0$. }
\medskip

{\it Proof of Lemma 6: } Clearly, the 
real and imaginary parts of $u$ are real analytic functions. 
Thus there is a Taylor expansion around $z=0$. Let 
\beq\label{taylor}
u(z,\ol{z}) 
= \sum_{|\alpha| =m} a_\alpha
z^{\alpha_1}\ol{z}^{\alpha_2} 
+ O\left(|z|^{m+1}\right),\ \ \ \alpha = (\alpha_1,\alpha_2)
\eeq
We claim that $a_\alpha =0$ except for the indices $\alpha = (m,0)$
and  $\alpha = (0,m)$. This statement is surely true if $m=1$. 
Henceforth we will assume that $m\geq 2$. 

We select $\alpha = (m_1,m_2)$. Let $\alpha ! = m_1! m_2!$. Then 
\beq
a_\alpha =  \frac{1}{\alpha!}
\left.\frac{\partial^{m_1+m_2}u}{ {\partial z}^{m_1}
{\partial\ol{z}}^{m_2}}\right|_{z=0} 
\eeq
Consider now the case that $m_1 \geq 1$ and $m_2 \geq 1$. Rewriting
\refeq{GLeq} as
\beq\label{GLeqC}
- 4 \frac{\partial^{}}{{\partial z}^{}}
\frac{\partial^{}}{{\partial\ol{z}}^{}} 
u(z,\ol{z})  = u(1-|u|^2),\ \ \ z\in \C
\eeq
and applying the $m_1-1$ st derivative in $z$ and the
 $m_2-1$ st derivative in $\ol{z}$ to \refeq{GLeqC}
 gives the identity
\beq\label{aa}
a_\alpha = -\frac{1}{4\alpha!} 
\left.\frac{\partial^{m_1-1}}{{\partial z}^{m_1-1}}
\frac{\partial^{m_2-1}}{{\partial\ol{z}}^{m_2-1}} \,
u(1-|u|^2)\right|_{z=0}  \, ,
\eeq
with $|\alpha | = m$, $\alpha = (m_1,m_2)$,
and $m_i \geq 1$ for $i = 1,2$.
Inserting \refeq{taylor} in \refeq{aa}, we find that $a_\alpha =0$ for 
all $\alpha$ satisfying $|\alpha| \leq m-1$. In particular, from
\refeq{aa},  $a_\alpha =0$ for 
$\alpha = (m_1,m_2)$ with $m_1 +m_2 = m$ and $m_1 \geq 1$, $m_2 \geq 1$.
Thus, 
\beq\label{uatori}
u(z,\ol{z}) =  a_{(m,0)} z^{m} + a_{(0,m)} \ol{z}^{m} 
+ O\left(|z|^{m+1}\right),
\eeq

Now we assume w.l.o.g. that $|a_{(m,0)}| \geq |a_{(0,m)}|$
and begin with the case 
 $|a_{(m,0)}| = |a_{(0,m)}|$. Then 
\beq
a_{(0,m)} = e^{i\theta_0}a_{(m,0) }.
\eeq
Thus, from \refeq{uatorirew}, we obtain now
\beq\label{uatorireduc}
u(z,\ol{z}) 
=  a_{(m,0)} \left( z^{m} +  e^{i\theta_0}
\ol{z}^m  \right) + O\left(|z|^{m+1}\right)
\eeq
which for $|z| = \eps$ gives
\beq\label{uatorieps}
u(z,\ol{z}) =  \left( 2 a_{(m,0)} e^{i\theta_0/2 }
 \cos{(m\theta + \theta_0/2)} + O(\eps )\right)\eps^m \, .
\eeq  
It follows now from \refeq{uatorieps} that the image of $|z|=\eps$
under $u$ has degree zero and has other zeros beside the one at the
origin. This is a contradiction, whence
 $|a_{(m,0)}| > |a_{(0,m)}|$. We now rewrite  \refeq{uatori} as 
\beq\label{uatorirew}
u(z,\ol{z}) =  a_{(m,0)} z^{m} \left[ 1 +\frac{a_{(0,m)}}{a_{(m,0)}} 
\frac{\ol{z}^m}{z^m} + O\left(|z|^{}\right)\right],
\eeq
It follows from \refeq{uatorirew} that, for $|z| = \eps >0$
sufficiently small, $u$ has an image curve whose winding number is
$m$. Thus, $m=d$. QED
\bigskip

As a consequence of Lemma 6, we have

\medskip
\noin
{\bf Corollary 1. -- } 
{\it Let $u$ be a solution of \refeq{GLeqD}, \refeq{Dbc},
\refeq{onenull}. Then 
\beq\label{vdef}
v(x) = v(z,\ol{z}) = u\bigl( z^{1/d} , \ol{z}^{1/d}\bigr) 
\in C^{1,\alpha}(B_R), 
\eeq
with $\alpha = 1/d$.}
\medskip

\medskip
\noin
{\bf Corollary 2. -- } 
{\it Let $u$ be a solution of \refeq{GLeqD}, \refeq{Dbc}, \refeq{onenull}.
Let $d^*$ be conjugate to $d$, i.e. $1/d^* =1- 1/d$. 
Then $v(x)$ defined in \refeq{vdef} satisfies
\beq\label{vEQ}
-\Delta v(x) = 4d^{-2} r^{-2/d^*} v(x) (1 -|v|^2)
\eeq
for $x\in \ B_R(0)\backslash\{0\}$, the singularity at $0$ being removable, 
together with
\beq\label{vBC}
v(x) = {x\over |x|};\ \ \ x\in \partial B_R(0)\, .
\eeq
}
\medskip

\medskip
\noin
{\bf Lemma 8. -- }{\it Let $u$ be a solution of \refeq{GLeqD},
\refeq{Dbc}, \refeq{onenull}. Let $v(x)$, defined by \refeq{vdef},
satisfy hypotheses \refeq{curlfree}, such that 
$v(x) = - \nabla \phi(x)$. Then $(1-|\nabla\phi|^2)/r^{2/d^*}$ 
is constant on the level curves of $\phi$. }
\medskip

{\it Proof of Lemma 8: } Substituting $v(x) = - \nabla \phi(x)$
into \refeq{vEQ} gives 
\beq\label{GLv}
  \nabla \big(-\Delta \phi (x)\big) = \, \kappa(x)\nabla\phi(x),
\eeq
where 
\beq\label{kappaphi}
\kappa(x)=\, 4d^{-2} r^{-2/d^*}(1-|\nabla\phi(x)|^2).
\eeq
Taking the curl of \refeq{GLv} we get
\beq\label{kappatimesphi}
\nabla \kappa(x) \times \nabla\phi(x) = 0.
\eeq
By \refeq{kappatimesphi}
and \refeq{onenull} it follows that $\kappa(x)$ is constant on 
the level curves of $\phi$. QED
\bigskip

With the help of the single zero hypothesis \refeq{onenull} we now
conclude that $\kappa(x)$ is a function of $\phi(x)$. 
This is proved as in section 3. As a consequence,
for solutions of \refeq{GLeqD} that are of the type \refeq{curlfree},
we can reduce \refeq{GLeqD} to a single scalar equation 
valid in $B_R(0)$.

\medskip
\noin
{\bf Lemma 9. -- }{\it Let $u$ satisfy the hypotheses of Theorem 2,
such that $v(x)$ given by \refeq{vdef} is well defined and
given by $v(x) = - \nabla \phi(x)$. Then there exists some 
increasing, differentiable function $g: \R\to\R^+$ such that
$\phi(x)$ solves the semi-linear elliptic PDE
\beq\label{PDEphiD}
-\Delta \phi(x) = g(\phi(x))\, , \ \ \ x\in B_R(0),
\eeq
with 
\beq\label{BCphiD}
\phi(x) = 0\, ,\ \ \ x\in \partial B_R(0).
\eeq
}
\medskip

{\it Proof of Lemma 9: } 
To prove \refeq{PDEphiD},
we repeat almost verbatim the proof of Lemma 3, up to 
\refeq{gdefine}, replacing $\R^2$ by $B_R(0)$ and
$|x|\to\infty$ by $|x|\to R$. Moreover, 
instead of $\lambda(x)$ now $\kappa(x)=\Theta(\phi(x))$, 
with $\kappa(x)$ given in \refeq{kappaphi}. Thus, 
\refeq{lambdaTheta} is replaced by
\beq\label{kappaTheta}
4d^{-2}r^{-2/d^*}(1- |\nabla \phi(x)|^2) = \Theta(\phi)(x) \, .
\eeq
With \refeq{kappaTheta} and \refeq{kappaphi} we can 
rewrite \refeq{GLv} in the form \refeq{nablaphieq}, 
with $x\in B_R(0)$, from where we arrive at \refeq{PDEphiD}.
Using now the fact that $v(x) = -\nabla\phi(x)$ determines $\phi(x)$
only up to an arbitrary additive constant, we see that we are
free to choose the constant so that 
\beq\label{infphi}
 \inf_{x\in  B_R(0)}\phi(x) =0 \, .
\eeq
Instead of \refeq{gdefine}, we now  choose
\beq\label{gdefineD}
g(t) = \int_{0}^t \Theta(s) ds\, . 
\eeq
Since $\Theta(s) \geq 0$, it follows that $g(t)$ is increasing and nonnegative.

We next prove \refeq{BCphiD}. Since $v(x) = -\nabla\phi(x)$, 
our boundary condition \refeq{vBC} becomes 
$\nabla\phi(x) = - x/|x|$ for $x\in \partial B_R(0)$,
the gradient understood in the limit  as $x\to \partial B_R(0)$
from inside. Thus, $|\nabla\phi(x)|=1$
for $x\in \partial B_R(0)$, which together with \refeq{kappaTheta}
gives $\Theta(\phi(x)) = 0$ for $x\in \partial B_R(0)$. We conclude
that $\phi(x)$ is constant on $\partial B_R(0)$. 
The maximum principle, applied to \refeq{PDEphiD}, gives 
that the infimum of $\phi$ is taken at  $\partial B_R(0)$, whence
with \refeq{infphi} we obtain \refeq{BCphiD}. QED
\bigskip

The radial symmetry of $\phi$ now follows by 
again comparing the isoperimetric inequality and Pohozaev's 
identity. In fact, such a proof for \refeq{PDEphiD}, 
\refeq{BCphiD} is given in \cite{Lions}. 
Having radial symmetry of $\phi$, it follows immediately
that $\phi^\prime(r) = -F_d(r^{1/d})$, with $F_d(r)$ satisfying the ODE 
stated in Theorem 2. This completes our proof of Theorem 2. 
\medskip

\subsection{An Open Question regarding Degree $d$ in $\R^2$.}

Beside Theorem 2, one rather would like to prove 
the analog of our Theorem 1 for degree $d$ solutions
in $\R^2$. Notice that Lemmata 6--9 and the Corollaries 1 and 2
hold with $B_R(0)$ replaced by  $\R^2$, and with 
\refeq{BCphiD} replaced by \refeq{phiasymp}. 
Thus, one is again lead to consider solutions of 
\refeq{phiPDE} and \refeq{phiasymp} that have only one 
critical point at $x_0$, but this time with 
\refeq{lambdaTheta} replaced by \refeq{kappaTheta}. 
Clearly, since $|\nabla\phi|$ is not a-priori constant on
the level curves of $\phi$, the type of  proof given in
subsection 3.1 does not apply. 

However, Lemma 4 holds with unchanged proof. Therefore, 
radial symmetry follows if one can prove  the analog 
of Lemma 5 for $d>1$. For this one needs to work directly
with \refeq{Hsurfint} as $t\to -\infty$. In fact, since
the weak form of Lemma 4 guarantees that
$\liminf_{t\to -\infty} H(\rho(t)) \geq 0$,  to prove
the analog of Lemma 5 it is sufficient to show that
$\limsup_{t\to -\infty} H(\rho(t)) \leq 0$. Thus, 
since for large negative $t$ the level curves are 
star shaped, we can use the Cauchy-Schwartz inequality to obtain 
\beq
H(\rho(t)) \leq 
{1 \over 2} \int_{\partial \Lambda_t}\langle x,\nu\rangle
 |\nabla \phi|^2 d\sigma 
\left(\int_{\partial\LL_t}
\frac{d\sigma}{\langle x,\nu\rangle} - 2\pi \right)\, .
\eeq
By Shafrir's asymptotics \refeq{moduasymp}, we infer that 
\beq\label{gradphiDasymp}
|\nabla\phi|^2 = 1 - \frac{d^2}{r^{2/d}} + o\left(\frac{1}{r^{2/d}}
\right)\, ,
\eeq
which is a weaker approach to 1 than for $d=1$, but strong enough to conclude 
\beq
{1 \over 2}\int_{\partial \Lambda_t}\langle x,\nu\rangle
 |\nabla \phi|^2 d\sigma
= |\LL_t| + o(|\LL_t|) = \pi \rho(t)^2 + o(\rho(t)^2)
\, .
\eeq
Therefore, the proof of the analog of Lemma 5 for $d>1$ is complete 
if one can show that
\beq\label{openprob}
\int_{\partial\LL_t} \frac{d\sigma}{\langle x,\nu\rangle} =
2\pi + o\left( \frac{1}{\rho(t)^2}\right)\, .
\eeq
This, however, is an open problem. 

\medskip
{\bf Acknowledgement. } We wish to thank H. Brezis for introducing
us to this problem and for his continued interest and encouragement. 
Thanks go to J.J. Aly for his comments on the
single zero condition, and to 
P. Mironescu for communicating his results to
us prior to publication, and for interesting discussions.
The work of S.C. was supported by NSF Grant \# DMS-9623079.
The work of M.K. was supported in parts by NSF Grant \# DMS-9623220.



\bigskip\bigskip
\begin{thebibliography}{[99]}

\bibitem{BandleBOOK}
Bandle, C.,
{\it Isoperimetric Inequalities and Applications},
 Pitman, Boston, (1980).

\bibitem{BethuelBrezisHelein}
Bethuel, F., Brezis, H., and H\'elein, F.,
{\it Ginzburg-Landau Vortices}, Birkh\"auser, Basel,
 (1994).

\bibitem{BrezisMerleRiviere}
Brezis, H., Merle, F., and Riviere, T.,     
{\it Quantization effects for $ -\Delta u = u (1 - |u|^2) \ in \ {\R}^2$},
Arch. Rat. Mech. Anal. {\bf 126}, 35--58 (1994).

\bibitem{brezisPC} H. Brezis, private communication. 

\bibitem{ckGAFA} 
Chanillo, S., and Kiessling, M.K.-H.,
{\it Conformally invariant systems of nonlinear PDE of Liouville type},
Geom. Funct. Anal. {\bf 5}, pp. 924--947 (1995).

\bibitem{ckCR}
Chanillo, S., and Kiessling, M.K.-H., 
Compt. Rend. Acad. Sci. Paris, {\bf t. 321}, S\'erie I., 1023--1026, (1995).

\bibitem{CohenNeuRosales} 
 Cohen, D.S., Neu, J.C., and Rosales, R.R., 
{\it Rotating spiral wave solutions of reaction-diffusion equations},
SIAM J. Appl. Math. {\bf 35}, 536--547, (1978).

\bibitem{GNN}
Gidas, B., Ni, W.M. and Nirenberg, L.,
{\it Symmetry and related properties via  the maximum principle},
Comm. Math. Phys. {\bf 68}, 209--243 (1979).

\bibitem{Greenberg} 
Greenberg, J.M., 
{\it Spiral waves for $\lambda-\omega$ systems},
SIAM J. Appl. Math. {\bf 39}, 301--309 (1980). 

\bibitem{Hagan}  
Hagan, P.S., 
{\it Spiral waves in reaction-diffusion equations},
SIAM J. Appl. Math. {\bf 42}, 762--786 (1982).

\bibitem{HerveHerve} 
Herv\'e, R.-M., and  Herv\'e, M., 
\'Etude qualitative des solutions 
r\'eelles d'une \'equation diff\'eren\-tiel\-le
li\'ee \`a l'\'equation de Ginzburg--Landau,
{\it Ann. Inst. Henri Poincar\'e} 11, 427-440, (1994).

\bibitem{JaffeTaubes}
Jaffe, A., and Taubes, C.,
{\it Vortices and Monopoles},
Birkh\"auser, Basel (1980).

\bibitem{KesavanPacella} 
Kesavan, C., and Pacella, F., 
{\it Symmetry of positive solutions of a quasilinear elliptic equation
via isoperimetric inequalities},
Appl. Anal. {\bf 54}, 27--37 (1994).

\bibitem{KopellHoward} 
Kopell, N., and Howard, L.N., 
{\it Target pattern and spiral solutions to reaction-diffusion equations with
more than one space dimension}, 
Adv. Appl. Math. {\bf 2}, 417--449 (1981).

\bibitem{LiebLoss}
Lieb, E.H., and Loss, M.,
{\it Symmetry of the Ginzburg--Landau Minimizer in a Disc},
Math. Res. Lett. {\bf 1}, pp. 701--715 (1994).


\bibitem{Lions} 
Lions, P.-L., 
{\it Two geometrical properties of solutions of semilinear problems},
Appl. Anal. {\bf 12}, 267--272 (1981).

\bibitem{MironescuCR}
Mironescu, P., 
{\it Les minimiseurs locaux pour l'\'equation de Ginzburg--Landau
sont \`a sym\'etrie radiale},
Compt. Rend. Acad. Sci. Paris, {\bf }, S\'erie I., 327--331
(1996); and  preprint (to appear.)


\bibitem{Shafrir}
Shafrir, I., 
{\it Remarks on solutions of $ -\Delta u = u (1 - |u|^2) \ in \ {\R}^2$},
Compt. Rend. Acad. Sci. Paris, {\bf t. 318}, S\'erie I., 327--331,
(1994).

\end{thebibliography}
\end{document}